\documentclass[11pt]{article}
\newcommand{\bbox}{~\hfill\rule{1.4mm}{3.2mm}}

\usepackage{multirow}

\def\Eq#1{Eq.\ (\ref{#1})}
\def\mod{\mathrel{\rm mod}}
\def\ydiv{\mathrel{|}}
\def\ndiv{\mathrel{\rlap{\kern-.05em\mbox{$\scriptstyle /$}}|}}

\newtheorem{defn}{Definition}
\newtheorem{lemma}{Lemma}
\newtheorem{theorem}{Theorem}

\title{\bf Divisibility tests with weighted digital sums}

\author{\bf Palash B. Pal\footnote{pbpal@theory.saha.ernet.in}\\
Saha Institute of Nuclear Physics\\
1/AF Bidhan-Nagar, Calcutta 700064, India}

\date{}

\begin{document}

\maketitle

\begin{abstract}

We propose a divisibility test for all integers which have 1, 3, 7 or
9 in their unit's place.  In particular, then, the test applies for
all prime divisors except 2 and 5.

\end{abstract}

\section{Introduction}
This article contains the statement and proof of an algorithm for
checking divisibility by all integers which have 1, 3, 7 or 9 in their
unit's place.  In particular, it contains divisibilty tests for all
primes except 2 and 5.  All algebraic symbols in this article
represent integers.

Our inspiration is obtained from a book by Kordemsky \cite{moscow},
where he discussed a test for checking divisibility by 7.  We present
it first for the sake of completeness, then proceed to generalize it
to other divisors.

Kordemsky's algorithm \cite{moscow} runs as follows.  Consider the
following decimal representation of a number $A$:
\begin{eqnarray}
A = \sum_{k}  a_k \times 10^k\,, \qquad \mbox{with $0\leq
  a_k \leq 9$} \quad \forall k \,. 
\label{A}
\end{eqnarray}
The sum has finite number of terms for any finite $A$.  In fact, if
there are $n+1$ digits in $A$, $a_k=0$ for all $k>n$.  We now have to
take $a_0$, i.e., the number in the unit's place, multiply by 5, add
$a_1$ to it, multiply the sum by 5, add $a_2$ to it, and continue like
this until we exhaust all digits of $A$.  At any point during the
operations, we can drop any multiple of 7.  If the final result is
divisible by 7, the original number $A$ is divisible by 7 as well.

Let us take an example to illustrate the procedure.  Consider
$A=3941$.  The steps in the algorithm can be represented in a tabular
form as follows:
\begin{eqnarray}
\begin{tabular}{ccccc}
\hline
\multirow{2}{1cm}{digit} &
\multirow{2}{2cm}{add result of last row} & \multirow{2}{13mm}{mod 7} &
\multirow{2}{1cm}{$\null \times 5$} & \multirow{2}{13mm}{mod 7} \\  
\\ 
\hline
1 &   & 1 & 5 & 5 \\ 
4 & 9 & 2 & 10 & 3 \\ 
9 & 12 & 5 & 25 & 4 \\ 
3 & 7 & 0 & 0 & 0 \\ 
\hline
\end{tabular}
\label{7table}
\end{eqnarray}
The final result is zero, and therefore $7\ydiv3941$, i.e., 3941 is
divisible by 7, as can be checked by direct division.

In \S\ref{proof}, we provide a generalization of this algorithm for
other divisors.  \S\ref{illvar} contains illustrations and comments.

\section{The proposed divisibility test and its proof}
\label{proof}
Consider any number which has 1, 3, 7 or 9 in its unit's place.  As
emphasized earlier, this includes all primes except 2 and 5, although
the algorithm described holds independently of whether the divisor is
a prime.  All such divisors have a multiple which has 9 in the unit's
place.  Thus, these divisors $d$ satisfy a relation of the form
\begin{eqnarray}
md = 10q - 1 
\label{dq}
\end{eqnarray}
for some integers $m$ and $q$.  We want to test whether $A$, given in
\Eq{A}, is divisible by $d$.  We have to follow an algorithm like the
one discussed by Kordemsky, except that now at each step we need to
multiply by $q$.  Formally, it is equivalent to Theorem \ref{thetest}
given below.  In order to build up to it, we need some definitions.

\begin{defn}
Given a number $A$ as in \Eq{A}, the numbers $\hat C_k$ are defined by:
\begin{eqnarray}
\hat C_0 &=& a_0 \,, \nonumber\\*
\hat C_k &=& q \hat C_{k-1} + a_k  \qquad \mbox{for $k\geq 1$}. 
\label{Chat}
\end{eqnarray}
\end{defn}
Obviously, if the number of digits in $A$ is $n+1$, we will obtain
$\hat C_k=\hat C_n$ for all $k>n$.  Thus, $\hat C_k$ will not change
if we increase $k$ indefinitely.  The terminal value of $\hat C_k$
will be denoted by $\hat C$, which can be called {\sl the $q$-weighted
digital sum of $A$}.  For a number $A$ having $n+1$
digits, the definition in \Eq{Chat} gives
\begin{eqnarray}
\hat C = q^n a_0 + q^{n-1} a_1 + \cdots + a_n \,.
\end{eqnarray}
\begin{defn}
Given a number $A$ as in \Eq{A}, the numbers $C_k$ are defined by:
\begin{eqnarray}
C_k &=& \hat C_k \mod d \qquad \forall k \,.
\label{C}
\end{eqnarray}
\end{defn}
Like its hatted relative, $C_k$ also has a terminal value, which will
be denoted by $C$, which can be called {\sl the $d$-modded digital sum
of $A$}.  In a practical situation, calculating $C_k$'s are much
easier than calculating $\hat C_k$'s, because $C_k$ involves
multiplication of smaller numbers in general.

\begin{defn}  
Suppose the number $A$, represented in \Eq{A}, is written
alternatively as
\begin{eqnarray}
A = \sum_{k=0}^n  a'_k \times 10^k\,, 
\label{Ab}
\end{eqnarray}
without any restriction on the $a'_k$'s except that all of them need
to be integers.  Then the ordered set $\{a'_k\}$ will be called a
``shuffled representation'' of the number $A$, and the process from
changing from one such representation to another will be called a
``shuffling''.
\end{defn}
Obviously, shuffled representations of a number include its regular
decimal representation, which has the additional restriction $0\leq
a'_k\leq 9$ for all $k$.  But $a'_k$ is not restricted to be a
single-digit positive number.  For example, consider the number 154,
for which $a_2=1$, $a_1=5$, $a_0=4$.  However, with the definition in
\Eq{Ab}, there are other options.  For example, we can choose
$a'_2=0$, $a'_1=15$, $a'_0=4$, or $a'_2=1$, $a'_1=3$, $a'_0=24$, or
even $a'_2=-1$, $a'_1=25$, $a'_0=4$.

\begin{lemma}\label{shuffinv}
The number $C$ is invariant under a shuffling.
\end{lemma}
\paragraph{Proof of Lemma \ref{shuffinv}~: } 
Any shuffling can be built up of succesive application of two kinds of
basic shufflings.  One of them is $S_r$, defined by
\begin{eqnarray}
a'_r = a_r-1 \,, \qquad a'_{r-1} = a_{r-1} + 10 \,,
\label{Sr}
\end{eqnarray}
and the other is $S_{\bar r}$, defined by
\begin{eqnarray}
a'_r = a_r+1 \,, \qquad a'_{r-1} = a_{r-1} - 10 \,,
\label{Srbar}
\end{eqnarray}
where in both definitions it is implied that $a'_k=a_k$ when $k\not=r$
or $k\not=r-1$.  Take $S_r$ first.  If we follow the procedure of
\Eq{C} using this shuffled set, suppose we obtain the numbers
$C_k^{(r)}$.  It is then enough to show that $C^{(r)}=C$.

Obviously, the shuffling $S_r$ does not disturb the numbers $a_k$ for
$k\leq r-2$, and therefore $C_k^{(r)}=C_k$ for $k\leq r-2$.  Using
\Eq{Sr}, we find that the next step gives
\begin{eqnarray}
C_{r-1}^{(r)} = \Big( C_{r-1} + 10 \Big) \mod d \,,
\end{eqnarray}
and, in the step after that, we get
\begin{eqnarray}
C_r^{(r)} &=& \Big( q(C_{r-1} + 10) + a_r -1 \Big) \mod d \nonumber\\
&=& \Big( \hat C_r + md \Big) \mod d \,,
\end{eqnarray}
using \Eq{dq} in the last step.  Since $md \mod d =0$, this can be
rewritten as
\begin{eqnarray}
C_r^{(r)} = \hat C_r \mod d = C_r \,.
\end{eqnarray}
For $k>r$, it does not make any difference whether we are using the
$a_k$'s or the shuffled set.   Thus $C^{(r)} = C$.  The proof is
similar for the basic shuffling $S_{\bar r}$, and so the lemma is
proved. \bbox 

\begin{lemma}\label{relprime}
The numbers $q$ and $d$, related through \Eq{dq}, are relatively prime.
\end{lemma}
\paragraph{Proof of Lemma \ref{relprime}~: } 
If $q=1$, there is nothing to prove.  For $q>1$, we will use the
result~\cite{book} that, for integers $a>b>0$,
\begin{eqnarray}
\gcd(a,b) = \gcd(a-b,b) \,,
\label{gcdrule}
\end{eqnarray}
which forms the basis of the Euclid algorithm for finding greatest
common divisors.  Repeated application of this rule gives
\begin{eqnarray}
\gcd(md,q) = \gcd(q-1,q) = \gcd(q,q-1) \,,
\end{eqnarray}
where in the last step, we have used the fact that the operation
$\gcd$ is commutative.  Using \Eq{gcdrule} again, we obtain
\begin{eqnarray}
\gcd(md,q) = \gcd(1,q-1) = 1 \,.
\end{eqnarray}
If there exists no common factor between $md$ and $q$, certainly there
does not exist any common factor between $d$ and $q$.  This completes
the proof of the lemma. \bbox

\begin{theorem}\label{thetest}
For divisors of the form given in \Eq{dq}, 
\begin{eqnarray}
C = 0 \Leftrightarrow d\ydiv A \,.
\label{9rule}
\end{eqnarray}
\end{theorem}
Note that for $d=9$, we have $q=1$ in \Eq{dq}, and in that case the
procedure of finding $C$ coincides with finding the digital sum of $A$
modulo 9, which is what constitutes the usual divisibility test for 9.
We now proceed for a general proof of the theorem.

\paragraph{Proof of Theorem \ref{thetest}~: }
We first prove that $d\ydiv A \Rightarrow C = 0$.  If $d\ydiv A$, we
can write $A=rd$ for some integer $r$.  We can take a shuffled
representation of $A$ in the form
\begin{eqnarray}
a'_0 = rd\,, \qquad a'_k=0 \quad \forall k\geq 1 \,.
\end{eqnarray}
Obviously then $C'_0 = rd \mod d=0$.  Since the higher $a'_k$'s are
zero, higher values of $C'_k$ remain zero, and so the terminal value
$C'=0$.  By Lemma \ref{shuffinv}, we then obtain $C=0$.

To prove the converse, suppose $d\ndiv A$.  Then $A=r \mod d$, with
$0<r<d$.  Then $B=A-r$ will be divisible by $d$.  A shuffled
representation of $B$ can be taken as
\begin{eqnarray}
b'_0 = a_0 - r \,, \qquad b'_k = a_k \quad \mbox{for $k\neq 0$}.
\end{eqnarray}
The $q$-weighted digital sum calculated from this representation will
be 
\begin{eqnarray}
\hat D = q^n (a_0-r) + q^{n-1} a_1 + \cdots + a_n \,.
\end{eqnarray}
Thus, 
\begin{eqnarray}
\hat C - \hat D = r q^n \,.
\end{eqnarray}
After modding out with $d$, we then obtain
\begin{eqnarray}
C - D = r q^n \mod d \,.
\end{eqnarray}
Since $d\ydiv B$ by definition, $D=0$ by the first part of the
proof.  Thus we have 
\begin{eqnarray}
C = r q^n \mod d \,.
\label{remainder}
\end{eqnarray}
By Lemma \ref{relprime}, $q$ and $d$ cannot have any common factor.
Also, $d\ndiv r$ since $0<r<d$.  Therefore we have proved that
$d\ndiv A \Rightarrow C\not=0$.  Inverting the logic, we thus obtain
$C=0 \Rightarrow d \ydiv A$, and the proof is complete.\bbox

\section{Illustrations and variations with prime divisors}
\label{illvar}
It is to be noted that the class of divisors defined in \Eq{dq}
contains all integers which have 1,3,7 or 9 in their unit's place.  In
particular, it contains all primes except 2 and 5.  Along with the
trivial tests for divisibility by 2 and 5 by checking the digit in the
unit's place, the present method then gives divisibility tests for all
primes.  We give some examples, with some important comments at the end.

We first give an illustration where $d$ is of the form $10q-1$.  Take
$A=16762$ and $d=29$.  We now have $q=3$.  Making a table as in
\Eq{7table}, we now obtain the following:
\begin{eqnarray}
\begin{tabular}{cccccc}
\hline
\multirow{2}{1cm}{digit} &
\multirow{2}{2cm}{add result of last row} & \multirow{2}{13mm}{mod 29} &
\multirow{2}{1cm}{$\null \times 3$} & \multirow{2}{13mm}{mod 29} \\  
\\ 
\hline
2 &   & 2 & 6 & 6 \\ 
6 & 12 & 12 & 36 & 7 \\ 
7 & 14 & 14 & 42 & 13 \\ 
6 & 19 & 19 & 57 & 28 \\ 
1 & 29 & 0 & 0 & 0 \\ 
\hline
\end{tabular}
\end{eqnarray}
Thus 16762 is divisible by 29.

Our next example will test the divisibility of 32893 by 7.  The
smallest positive value of $m$ satisfying \Eq{dq} is $m=7$, which
gives $q=5$, as in Kordemsky's original algorithm \cite{moscow}.  The
7-modded digital sum of 32893 is now computed in the following table.
\begin{eqnarray}
\begin{tabular}{cccccc}
\hline
\multirow{2}{1cm}{digit} &
\multirow{2}{2cm}{add result of last row} & \multirow{2}{13mm}{mod 7} &
\multirow{2}{1cm}{$\null \times 5$} & \multirow{2}{13mm}{mod 7} \\  
\\ 
\hline
3 &   & 3 & 15 & 1 \\ 
9 & 10 & 3 & 15 & 1 \\ 
8 & 9 & 2 & 10 & 3 \\ 
2 & 5 & 5 & 25 & 4 \\ 
3 & 7 & 0 & 0 & 0 \\ 
\hline
\end{tabular}
\end{eqnarray}
This shows that 32893 is divisible by 7.  However, it should be
commented that it is not essential that we take the smallest positive
integers satisfying \Eq{dq}.  For example, if we take $m=-3$ in
\Eq{dq}, we obtain $q=-2$.  We could have easily used this value of
$q$ to obtain the modded digital sum, as shown below.
\begin{eqnarray}
\begin{array}{cccccc}
\hline
\multirow{2}{1cm}{digit} &
\multirow{2}{2cm}{add result of last row} & \multirow{2}{13mm}{mod 7} &
\multirow{2}{15mm}{$\null \times (-2)$} & \multirow{2}{13mm}{mod 7} \\  
\\ 
\hline
3 &   & 3 & -6 & -6 \\ 
9 & 3 & 3 & -6 & -6 \\ 
8 & 2 & 2 & -4 & -4 \\ 
2 & -2 & -2 & 4 & 4 \\ 
3 & 7 & 0 & 0 & 0 \\ 
\hline
\end{array}
\end{eqnarray}
Clearly, the result is the same.

This freedom of choosing $m$ and $q$ is particularly helpful for some
larger primes.  For example, if $d=11$, smallest positive integers
satisfying \Eq{dq} are given by $m=9$, $q=10$.  It is much easier to
work with $m=-1$, $q=-1$, which is equivalent to the usual
divisibility test for 11.  For $d=17$, smallest positive integers
satisfying \Eq{dq} are given by $m=7$, $q=12$.  Finding the modded
digital sum, one therefore needs to multiply by 12 at each stage.  It
is much easier to take the solutions $m=-3$, $q=-5$ instead, where the
computation of the modded digital sum will involve multiplication by
$-5$.

Finally, a word of caution.  The algorithm described here can test
whether a number is divisible by a certain divisor.  If it is not, the
algorithm does not provide a quick answer for what the remainder might
be.  As \Eq{remainder} shows, if $d\ndiv A$, the relation between the
modded digital sum $C$ and the remainder $r$ depends on the number of
digits in $A$, unless of course $q=1$.



\begin{thebibliography}{[99]}

\bibitem{moscow} Boris A. Kordemsky, {\sl Moscow puzzles : 359
mathematical recreations}; edited and with an introduction by Martin
Gardner.  Dover Publications, New York, 1992.  Problem 320 contains
the algorithm mentioned here.


\bibitem{book}  See any textbook on number theory, e.g., Song Y. Yan,
  {\sl Number theory for computing}; Springer, 2nd edition, 2002.

\end{thebibliography}
\end{document}